\documentclass[a4paper,11pt]{amsart}
\usepackage{amsfonts}
\usepackage{amsmath}
\usepackage{amssymb}
\usepackage{amsthm}
\usepackage{young}
\usepackage[all]{xypic}

\date{\today}
\newtheorem{teor}{Theorem}[section]
\newtheorem{coro}[teor]{Corollary}
\newtheorem{conj}[teor]{Conjecture}

\newtheorem{lem}[teor]{Lemma}
\newtheorem{prop}[teor]{Proposition}
\newtheorem{defi}[teor]{Definition}
\newtheorem{example}{Example}
\newtheorem{remark}{Remark}

\newcommand{\R}{\mathbb R}

\newcommand{\Q}{\mathbb Q}
\newcommand{\C}{\mathbb C}

\newcommand{\cqd}{\hfill$\Box$}

\def\diso{\lower.4ex\hbox{$\downarrow$}\raise.4ex\hbox{\mbox{\scriptsize $\wr$}}}

\def\genk#1{\langle\, {#1}
\,\rangle_{\lower.4ex\hbox{{\scriptsize$K$}}}}
\def\genq#1{\langle\, {#1}
\,\rangle_{\lower.4ex\hbox{{\scriptsize$\Q$}}}}
\def\genr#1{\langle\, {#1}
\,\rangle_{\lower.4ex\hbox{{\scriptsize$\R$}}}}

\def\iso{\,\lower .6ex\hbox{$\stackrel{\lra}{\mbox{\tiny $\sim\,$}}$}\,}

\def\lg{l\raise.6ex\hbox to.2em{\hss.\hss}l}
\def\lra{\longrightarrow}

\def\orb{\hbox to  .3em{$\backslash$}\backslash}

\newenvironment{pr}{\noindent{\llap{\bf\theprobl .}}\ \ }
                   {\stepcounter{probl}\vspace{.cm}}
\newcommand{\p}{\begin{pr}}
\newcommand{\fp}{\end{pr}}
\newcounter{probl}
\stepcounter{probl}

\newcounter{cs}
\stepcounter{cs}
\newcommand{\casos}{\begin{itemize}}
\newcommand{\fcasos}{\end{itemize}\setcounter{cs}{1}}

\newfont{\tit}{cmr12 scaled \magstep3}

\title[Counting cyclic coverings of $\mathbb{P}^{1}$.]{ Enumerative geometry  of the curves defined by $y^{n}=f(x)$.}

\author[Cristina Mart{\'\i}nez]{Cristina Mart{\'\i}nez}
\author[Alberto Besana]{Alberto Besana}

\address{
Departament de Matematiques, Edifici C, Facultad de Ciencies, Universitat Aut\`onoma de Barcelona, 08193 Bellaterra, Barcelona}
 \email{cmartine@mat.uab.cat}
 \email{abesana@mat.uab.cat}
   \address{Dipartimento di Matematica "Guido Castelnuovo"
Sapienza Universit\`a di Roma
P.le Aldo Moro, 5 - 00185 Roma}
\email{martinez@mat.uniroma1.it}

\begin{document}
\maketitle

\begin{abstract}
We study plane algebraic curves defined over a field $k$ of arbitrary characteristic as coverings of the the projective line $\mathbb{P}^{1}(k)$ 
and  the problem of enumerating branched coverings of the projective line by using combinatorial methods.
\end{abstract}

{\small \tableofcontents }

\section{Introduction}
There are strong analogies between plane curves and coverings of the projective line defined over a  field $k$ of arbitrary characteristic. In the present paper we study their connections and their relations with the combinatorics of Hurwitz numbers.
In particular, 
we study curves of the form 
\begin{equation}\label{equ1}
y^{d}=(x-a_{1})^{m_{1}}\ldots (x-a_{r})^{m_{r}},
\end{equation} 
for given $d,r$ and  $m_{i}$, integer numbers. It is easy to see that these curves correspond to coverings of $\mathbb{P}^{1}$ with Galois group $\mathbb{Z}_{d}$ acting by multiplication with a $d-$root of unity on the coordinate $y$ and with ramification at the points $a_{i}$. The data defining such a covering are encoded by a partition of length $d$. 
If $d,r$ and $m_{i}$ are coprime numbers,  the corresponding field extension $k(x)\hookrightarrow k(x,y)$ is a Kummer extension of the rational function field $k(x)$.


The Galois group of the plane curve  $C_{f}$ with affine model defined by the equation (\ref{equ1}) is the Galois group of the polynomial 

\noindent $f(x)=(x-a_{1})^{m_{1}}\ldots (x-a_{r})^{m_{r}},$ that is, the automorphism group $Aut\,(k(R_{f})/k))$, where $R_{f}$ denotes the set of branch points of the associated covering map $\pi: C_{f}\rightarrow \mathbb{P}^{1}$.
The Galois group of the curve $C_{f}$ is a quotient of the automorphism group $Aut(C_{f})$ of the curve and if it contains a cyclic subgroup $ \mathbb{Z}_{p}$, where $p$ is a prime number, such that the quotient curve $C_{f}/\mathbb{Z}_{p}$ has genus 0, then the curve is called a cyclic $p-$ gonal curve. If in addition $\mathbb{Z}_{p}$ is normal in $Aut(C_{f})$, then $C_{f}$ is called a normal cyclic $p-$gonal curve. In this case, the reduced automorphism group $\overline{Aut}(C_{f}):=Aut(C_{f})/\mathbb{Z}_{p}$  is isomorphic to a finite subgroup of $PGL_{2}(k)$.


We  study curves with Galois group $S_{n}$ and their invariant fields under the action of finite subgroups of $S_{n}$. 
In particular, we consider the locus of curves $X$ with reduced automorphism group isomorphic to the dihedral group $D_{n}$ such that $X/D_{n}\cong \mathbb{P}^{1}$. 
 
Let $\mathcal{C}_{d,m}$ the variety parametrizing the curves $C_{f,d}$, that is, the parameter space of coefficients of the equations of the form (\ref{equ1}). This is a Zariski open set in $\mathbb{A}_{k}$ corresponding to the complement $V(D)$, where $D$ is a suitable discriminant and itself an algebraic variety with coordinate ring $k[x_{1},\ldots, x_{d}]_{D}$. All the curves corresponding to points in $\mathcal{C}_{d,m}$ have the same genus $g$. The moduli space $H_{g,d}$ of pairs $(C_{f,d}, \pi: C\rightarrow \mathbb{P}^{1}) $ 
 is a Hurwitz space.

In Theorem \ref{classTh} we give a complete classification of all cyclic coverings of $\mathbb{P}^{1}$ over a field $k$ of characteristic $p\geq 0$, by  genus $g$ and degree $d$ curves with prescribed ramification over $\infty$ given by a partition of $d$.

In section \ref{sec4}, we study the enumerative problem of counting degree $d$ coverings of $\mathbb{P}^{1}$ by distinguishing on the number of ramification points. The enumeration of coverings of the complex projective line with profile $\mu$ over $\infty$ and simple ramification over a fixed set of finite points is done by direct calculation in the Gromov-Witten theory of $\mathbb{P}^{1}$. These numbers are known as Hurwitz numbers and arise as intersections in $\overline{M}_{g,n}(\mathbb{P}^{1})$.  
\subsubsection*{Conventions}
For $d$ a positive integer,  $\alpha=(\alpha_{1},\ldots, \alpha_{m})$ is a partition of $d$ into $m$ parts if the $\alpha_{i}$  are positive and non-decreasing. We set $l(\alpha)=m$ for the length of $\alpha$, that is the number of cycles in $\alpha$, and $l_{i}$ for the length of $\alpha_{i}$. The notation $(a_{1},\ldots, a_{k})$ stands for a permutation in $S_{d}$ that sends $a_{i}$ to $a_{i+1}$.
For us, scheme means separated scheme of finite type over an algebraically closed field $k$. 
A curve is an integral scheme of dimension 1, proper over $k$. 
We write $PGL(2,k)=GL(2,k)/k^{*}$, and elements of $PGL(2,k)$ will be represented by equivalence clases of matrices 
$\left(\begin{array}{ll} a &  b \\
c & d \end{array}\right)$, with $ad-bc\neq 0$.

We will denote the greatest common divisor of two integers $a$ and $b$ as $(a,b)$.

\section{Polynomial invariants under the action of a finite group}
Let $k$ be an algebraically closed field of characteristic $p\geq 0$ . 
Let V be a finite dimensional $k$-vector space equipped with a linear action, that is, $G$ acts via a representation $G \rightarrow GL(V )$. As an algebraic 
variety, $V$ is an affine space and its Picard group is  trivial. We denote by $k[V]$ the polynomial algebra over $k$ on $V$, and let $k[V,X]$ be the result of adjoining a formal variable $X$ to $k[V]$. We grade $k[V]$ by demanding the elements of $V$ to have degree one. So if $z_{1}\ldots,z_{n}$ is a basis for $V$ then $k[V]=k[z_{1},\ldots,z_{n}]$, where $z_{1},\ldots, z_{n}$ may also be thought of as formal variables of degree one.

The action of $G$ induces a natural action on the polynomial ring $k[V]\cong Sym(V^{*})$. The coordinate ring of invariant polynomials $k[V]^{G}$, 
is finitely generated as an algebra, for some homogeneous polynomials called $G-$ invariants. The locus $V(f_{1},\ldots,f_{r})$ defined by the invariant polynomials is an algebraic variety $X$ with coordinate ring $k[z_{1},\ldots,z_{n}]{G}$. The function field of $X$ is defined as the quotient field $k(z_{1},\ldots, z_{n})/(f_{1}\ldots, f_{r})$, where $k(z_{1},\ldots, z_{n})$ is the function field of the projective space $\mathbb{P}^{n}(k)$, and $(f_{1},\ldots, f_{r})$ is the ideal generated by the polynomials $f_{1},\ldots, f_{r}$.

 When $G$ has a polynomial ring of invariants, we define the Jacobian determinant $J=J(f_{1},\ldots,f_{n})=det(\frac{\partial\, f_{i}}{\partial\, z_{j}})$. This polynomial is nonzero and well-defined up to a nonzero element of $\mathbb{C}$ depending on the choice of basic invariants of a basis $\{z_{j}\}$ of $V^{*}$. 

When does $G$ have a polynomial ring of invariants?
Serre showed that in arbitrary characteristic, every finite subgroup of $GL(V)$ with a polynomial ring of invariants must be generated by reflections (see \cite{Se}). The converse may fail when the characteristic of the field divides the order of $G$.

The ring of polynomials in $n$ variables with complex coefficients admits a natural action of the orthogonal group $SO(n)$. We can also study the action of finite subgroups $G$ of $SO(n)$ and give generators for the spaces $\C[x_{0},\ldots, x_{n}]^{G}_{j}$ of homogeneous $G-$invariant polynomials of degree $j$. We can even compute their dimension by considering the Poincar\'e series
$$p(t):=\sum_{j=0}^{\infty} dim\,\, \C[x_{0}\ldots, x_{n}]^{G}_{j}\cdot t^{j}.$$

It can be written as $$p(t)=\frac{1}{|G|}\sum
\frac{n_{g}}{det(g-1\cdot t)},$$ where the sum runs over all the
conjugacy classes of $G$ and $n_{g}$ denotes the number of their elements.

We define the polynomials $f(x)=x^{n}+a_{1}x^{n-1}+\ldots+a_{n}$, with
$a_{i}\in \mathbb{Q}$, and
$f(x,t)=x^{n}+a_{1}x^{n-1}+\ldots+a_{n}-t$. Then, if $f$ is a
separable polynomial, then the Galois group of $f(x,t)$ over $k(t)$ is a regular extension with Galois group $S_{n}$. 
 
\begin{example}

\rm{Let $G=S_{n}$ acting on $\mathbb{Q}(x_{1},\ldots, x_{n})$. 
Observe that $\mathbb{Q}(x_{1}\ldots, x_{n})$ is the function field of an $(n-1)-$di\-men\-sio\-nal projective space $\mathbb{P}^{n-1}(\mathbb{Q})$ over $\mathbb{Q}$. Suppose that $z_{1}\ldots, z_{n}$ are the roots of $f$ in a splitting field of $f$ over $\mathbb{Q}$.  Each coefficient $a_{i}$ of $x^{i}$ in $f$ is symmetric in $z_{1}\ldots,z_{n}$, thus by the theorem on symmetric functions, we can write $a_{i}$ as a symmetric polynomial in $z_{1}\ldots, z_{n}$ with rational coefficients.
On the other side, for a permutation $\sigma\in S_{n}$, set $E_{\sigma}=x_{1}z_{(\sigma(1))}+\ldots +x_{n}z_{(\sigma(n))}$ in $\mathbb{Q}(x_{1}\ldots, x_{n})$ and $f(x)=\prod_{\sigma}(x-E_{\sigma})$, where $\sigma$ runs through all permutations in $S_{n}$.

\begin{teor} \label{Th1} (Serre)
The field $E$ of $S_{n}-$ invariants is $\mathbb{Q}(t_{1},\ldots,t_{n})$, where $t_{i}$ is the $i$th symmetric polynomial in $x_{1},\ldots, x_{n}$, and $\mathbb{Q}(x_{1}\ldots, x_{n})$ has Galois group $S_{n}$ over 
$E$: it is the splitting field of the polynomial
$$f(x)=x^{n}-t_{1}x^{n-1}+t_{2}x^{n-2}+\ldots+(-1)^{n}t_{n}.$$
\end{teor}

}
\end{example}

\begin{example}{Dihedral invariants}
\rm{ A dihedral group is the group of symmetries of a regular polygon, including both rotations and reflections. The Dihedral group $D_{s}$ is generated by a rotation $\tau$ of order $s$, and a reflection $\sigma$ of order 2, such that $\sigma \tau \sigma= \tau^{-1}$. In geometric terms, in the mirror a rotation looks like an inverse rotation. The action of the Dihedral group $D_{s}$ on $\mathbb{C}(x_{1},\ldots, x_{s})$ is given by 

 $$\tau: x_{j}\mapsto \epsilon^{i}x_{j}, \ \ for\ \  i=1,\ldots , [\frac{s}{2}].$$

$$\sigma:  x_{i}\mapsto x_{s-i}, \ \ for\ \ i=0,\ldots, [\frac{s}{2}], $$

where $\epsilon$ is an $s$ primitive root of unity.

If $s=1$, then the above actions are trivial. 
If $s=2$, then $\tau(x_{1})=-x_{1}$, $\tau_{1}(x_{2})=x_{2}, \tau_{2}=Id$, and the action is not dihedral but cyclic on the first factor.

We need to find the invariant polynomials in the coordinates $(x_{1},\ldots,x_{s})$ by the action of the Dihedral group. Let $s>2$, then the elements 
$$a_{i}(x_{1},\ldots, x_{s}):=x_{1}^{s-i}x_{i}+x^{s-i}_{s-1}x_{s-i}, $$
$$a_{s-i}(x_{1},\ldots, x_{s}):=x_{1}^{i}x_{s-i}+x^{i}_{s-1}x_{i}$$ for $1\leq i\leq s,$ are invariant  polynomials under the action of the group $D_{s}$ defined above. The elements $x_{i}$ are called the  dihedral invariants of $D_{s}$}.
\end{example}

\section
{ Cyclic coverings of $\mathbb{P}^{1}$ with prescribed ramification}
Let $k(x)$ be the function field of the projective line and consider a finite Galois extension $E$ of $k(x)$ with group $G$ which is regular, i.e., $\bar{k}\cap E=k$. \begin{lem}Let $G=Gal(E/k(x))$ be the Galois group of the extension. The inclusion $k(x) \hookrightarrow E$ corresponds to a (ramified) Galois covering $C\rightarrow \mathbb{P}^{1}$ defined over $k$ with Galois group $G$.  
\end{lem} \label{lema1}
{\it Proof.}
Geometrically, $E$ can be viewed as the function field $k(C)$ of a smooth projective curve $C$ which is absolutely irreducible over $k$, i.e., $y$ satisfies an algebraic equation over $k[x,t]$
$$\{(x,t)\in k^{2}|\,\, f(x,t)=0\},$$ where $f(x,t)=\sum_{i=1}^{m}\sum_{j=1}^{n}x^{i}y^{j}=\sum_{j=1}^{n}a_{j}(x)y^{j}=0$ is an irreducible polynomial in $x$ and $y$. We assume that not all $a_{ij}$ vanish and that $a_{n}(x)=1$ (which can be arranged by a change of variables). Thus $n$ is the degree of the polynomial in $y$. Since $k$ is algebraically closed, at a generic point $x$ there are $n$ roots $y^{(k)}$, $k=1,\ldots, n$ which implies that the algebraic curve defines an $n-$sheeted ramified covering of the $x-$plane given by projecting  over the $x$-axis. If the number of distinct roots $y^{(k)}$ is lower than the degree $n$, this means there are roots that occur with multiplicity greater or equal than 2. These are called branch points which belong to several sheets of the covering.

\cqd

\begin{conj} Every finite group $G $ occurs as the Galois group of such a covering.
\end{conj}

\begin{defi} Two coverings $\pi:C\rightarrow \mathbf{P}^{1}$ and $\pi': C'\rightarrow \mathbf{P}^{1}$ are isomorphic if there exists an isomorphism of curves $\phi: C\rightarrow C'$ satisfying $\pi'\circ \phi=\pi$.
\end{defi}

Define the special locus of a map $f:X\rightarrow \mathbf{P}^{1}$ (where $X$ is a nodal curve) as a connected component of the locus in $X$ where $f$ is not \'etale. Then a special locus is  a 
singular point on $X$ that is an $m-$fold branched point (analytically, the map looks locally like $x
\rightarrow x^{m}$, $m>1$), or a node of $X$, where the two branchings of the node are branch points of order $m_{1}, m_{2}$, or a one dimensional scheme of arithmetic genus $g$, attached to $s$ branches of the remainder of the curve that are $c_{j}-$fold branch points $(1\leq j \leq s)$. The form of the locus, along with the numerical data, will be called the type. Following \cite{GV}, to each special locus, associate a ramification number as follows:
\begin{itemize}
\item $m-1$
\item $m_{1}+m_{2}$
\item $2g-2+2s+\sum_{j=1}^{s}(c_{j}-1)$

\end{itemize}
The total ramification above a point of $\mathbb{P}^{1}$ is the sum of the ramification numbers of the special loci mapping to that point.
For any map $f$ from a nodal curve to a non-singular curve, the ramification number defines a divisor on the target:
$$\sum_{L}r_{L}f(L),$$where $L$ runs through the special loci and $r_{L}$ is the ramification index.

\begin{defi} Given a covering $\pi: C\rightarrow \mathbb{P}^{1}$ of degree $d$, the profile of $\pi$ over a point $q\in \mathbb{P}^{1}$ is the partition $\eta$ of $d$ obtained by the multiplicities of $\pi^{-1}(q)$.
\end{defi}

\begin{defi} Two ramified coverings $(C_{f}; \pi_{f}) $ and $(C_{g}; \pi_{g})$ are called topologically equivalent if there exists a homeomorphism $h: C_{f}\rightarrow C_{g}$ making the following diagram commutative:

$$ \xymatrix{ C_{f}  \ar[rr]^{h}  \ar[rd]^{\pi_{f}} &  & C_{g} \ar[ld] _{\pi_{g}}    \\
&  \mathbb{P}^{1} & } $$

In particular, the ramification points of the coverings coincide, as do the genera of the covering curves.
\end{defi}



Given a polynomial $f$  in $k[x]$ of degree $m$ with roots $\beta_{1}\ldots, \beta_{n}$ repeated according to the multiplicity in the splitting field $L$ of the extension of $f(x)$ over $k$,
and a positive integer $d$, let $C_{f,d}$ be the smooth projective curve over $k$ with affine model 
\begin{equation} \label{eq1}
y^{d}=f(x).
\end{equation}

We denote by $\xi_{d}$ a primitive $d-$th root in $\overline{k}$. Over $k(\xi_{q})$, we have a natural action of the $d-$th roots of unity on $C_{f,d}$, namely $(x,y)\rightarrow (x, \xi_{q}y)$. In par\-ti\-cular, we have that the Galois group of the covering is $Gal(k(C_{f,d}/k(x))\cong \mathbb{Z}_{d}$. This is why we call the extension a cyclic covering.
It is a Galois covering of the projective line that ramifies exactly at the places $x=\beta_{i}$, and the corresponding ramification indices are defined by $$e_{i}=\frac{n}{(n,d_{i})},$$
with $d_{i}$ the corresponding multiplicity of $\beta_{i}$ in $f$. There are ramification points with different ramification behavior. 

If $n \equiv 0\  (d)$, where $n:=\sum_{i=1}^{s}d_{i}$, then the place at $\infty$ does not ramify at the above extension. The only places of $k(x)$ that are ramified are the places $P_{i}$ that correspond to the points $x=\beta_{i}$. If the curve ramify at $\infty$, then $d_{i}\geq 2$ and the monodromy is given by a partition $\alpha$ of $d$. If the multiplicity $d_{j}=1$, then the point simply ramifies and the monodromy above the point is induced by a simple transposition. 
 

By the Riemann-Hurwitz formula, it follows that the function field $F$ has genus:
$$g=\frac{(n-1)(s-2)}{2}.$$

We denote by $R_{f}$ the set of roots of $f$ in $\overline{k}$. The function field of the curve is $F=k(x,y)$ where $y$ satisfies the algebraic equation ($\ref{eq1}$) over the algebra $k[x]$, that is, $k(C_{f})=k(x,y)$. Observe that $Gal(F/k(x)) \leqslant Aut(F) $.

We can consider the quotient surface $C/G$ for any finite subgroup $G$ of the automorphism group $Aut(F)$ of the curve  $C_{f,d}$. As we have seen, $C$ admits an automorphism $\tau$ of order $d$ such that $C/<\tau>$ is isomorphic to $\mathbb{P}^{1}$. The quotient surface is obtained via uniformizing a neighborhood of 0 by $y\rightarrow y^{d}$, this means the surface has at least an orbifold point of order $d$. The uniformization induces naturally an orbifold structure on the hyperplane class
bundle, such that the cyclic group $Z_{d}$  acts trivially on the corresponding bundle. The resulting orbifold bundle is denoted by $\mathcal{O}^{unif} (1)$.


\begin{defi} The Galois group of the curve $C_{f,d}$ is defined as the Galois group $Gal(f(x))$ of the polynomial $f(x)$, that is, the automorphism group $Aut(k(R_{f}/k))$. \end{defi}

\begin{defi} The discriminant of  the polynomial $f$ is $\triangle=\delta^{2}$, where 

$\delta=\Pi_{1\leq i<j\leq n} (\beta_{j}-\beta_{i})$. \end{defi}
If $f$ has a repeated root, then $\delta=0$, 
and $f$ is a separable polynomial if and only if $\delta\neq 0$.


\begin{lem}
If the base field $k$ is of characteristic 0 and $f(x)$ is an irreducible polynomial then $Gal(f(x))\cong S_{d}$.
\end{lem}
{\it Proof.} Just observe that being $f(x)$ an irreducible polynomial in a unique factorization domain $k[x]$, where $k$ is of characteristic 0, it is a separable polynomial. Thus there are $d$ different roots, where $d$ is the degree of the polynomial and the symmetric group $S_{d}$ acts by permuting them. In this case we say that the ramification is simple.\cqd

\begin{coro}

If $f(x)$ is a degree $d$ polynomial its Galois group $Gal(f(x))$ is a subgroup of the permutation group $S_{d}$ of $d$ elements.
\end{coro}

\begin{remark} Alternatively $C_{f,d}$ may be seen as an unramified Galois covering of a Riemann surface.  To describe the associated Riemann surface, one has to be able to identify the branching structure of the curve at the branch points, that is, one has to specify which sheets of the covering are connected in which way at a given branch point. This is equivalent to identifying the monodromy of the surface. Moreover every Riemann surface arises as a quotient of one of the simply connected domains $\mathbb{H}, \mathbb{C}$ and $\mathbb{P}^{1}$ by a discrete subgroup of the group of its automorphisms. These discrete subgroups are the fundamental groups of the corresponding underlying Riemann surface (see \cite{BB}). 
The cyclic coverings studied here  have genus greater or equal to 2 and are uniformized by the hyperbolic plane. Only rational curves have universal covering the projective line and only elliptic curves have universal covering the complex plane.

\end{remark}

Let $F_{0}$ be the fixed field $F^{\mathbb{Z}_{d}}$ by the action of the cyclic group $\mathbb{Z}_{d}$, then $Gal(F/F_{0})$ is the Galois group of the curve $C_{f}$ that is the Galois group of the polynomial $f(x)$.

Let  $Gal(F/k(x))=<\sigma>$ with $\sigma$ a generator of the Galois group, if $\tau\notin Gal (F/k(x))$, $\tau$ is said an extra automorphism. There is an exact sequence:
$$1\rightarrow \mathbb{Z}_{d}\rightarrow G\stackrel{\pi}{\rightarrow}G_{0}\rightarrow 1,$$
where $G_{0}=Aut(F)/\mathbb{Z}_{d}$ and $G=Aut(F)$. Moreover if the extension splits then $G_{0}\cong Gal(F/F_{0})$ and $Aut(F)  \cong \mathbb{Z}_{d}\times G_{0}$.

Let $b=div(f(x))_{0}$ be the root divisor of the polynomial $y^{d}=f(x)$ in $k(x)$.
The ramifications are determined by the profile of the covering over the branch points. Any branch point is induced by a permutation in $S_{d}$. In particular, if a point is symply ramified its monodromy is determined by a symple transposition. If we vary a branch point of the curve $C$ in $\mathbb{P}^{1}$, we obtain a one dimensional Hurwitz space parameterizing  such coverings. Each conjugacy class in $S_{d}$ determines a divisor class in  the Hurwitz space of all degree $d$ and fixed genus $g$ connected coverings of $\mathbb{P}^{1}$.

\begin{defi}A ramification type is realizable if the Galois group $Gal(F/F_{0})$ 
is a normal subgroup in the whole automorphism group $G=Aut(F)$.
\end{defi}
Observe that any normal finite subgroup $G_{0}$ of $Gal(F/F_{0})$ determines a ramification type.



\begin{lem} \label{lema}Any partition $\lambda=(\lambda_{1}, \ldots, \lambda_{m})$ of $d$ into $m$ parts corresponds to a degree $d$ branched covering of $\mathbb{P}^{1}$ with monodromy above $\infty$ given by $\alpha$, and $r=d+m+2\,(g-1)$ other simple branch points and no other branching. 
\end{lem}
{\it Proof.} 
For each partition $\lambda=(\lambda_{1}\geq\ldots \geq \lambda_{m} )$ in  $\mathcal{P}(d)$, consider  a configuration of points $\{p_{1},\ldots, p_{m}\}$ on the $x-$axis 
with coordinates
$$\{(x_{1},0), \ldots, (x_{m},0)\}\subseteq k^{*}\times \{0\} \subset k^{*}\times k . $$

To this configuration of points corresponds a unique polynomial

 $f(x)=(x-x_{1})^{\lambda_{1}}\ldots (x-x_{m})^{\lambda_{m}} \in k[x]$ which defines a covering of $\mathbb{P}^{1}(k)$. 
The divisor
$D_{\lambda}=\sum_{i=1}^{m}\lambda_{i}p_{i}$ 
corresponds to a ramification type defining the profile of the covering at $\infty$.
Every permutation $\sigma\in S_{d}$ defines an automorphism of the covering acting by permuting the places corresponding to the points $p_{i}, \ \  i=1,\ldots ,m$. 
In particular, permutations in the same conjugacy class have the same cycle structure and thus give the same ramification type.
\cqd

\begin{remark}
 Observe that if the base field $k$ is of characteristic different from 0, then not any configuration of points as in \ref{lema} gives rise to a polynomial in $k[x]$.
\end{remark}

\begin{defi} A set of integers $R$ $mod\ n$ is said a set of roots, if it is the set of roots of some polynomial, that is, if it corresponds to $R_{f}$ for some polynomial $f\in \mathbb{Z}[x]$.
\end{defi}

Let $q=p^{n}$ and consider the Galois extension $\mathbb{F}_{q}/\mathbb{F}_{p}$ with Galois group the cyclic group of order $n$. 
According to the Chinese remainder theorem, finding and counting sets of roots $mod \ n$ reduces to compute roots modulo a prime power (see \cite{Mau}). Indeed there is a functorial correspondence between polynomials in $\mathbb{Z}[x]$ modulo a prime and root sets. 

On the other hand the set of roots of a polynomial over $\mathbb{Z}$ coincides with the set of roots of a polynomial over $\mathbb{Q}$, that is, every rational root of a polynomial in $\mathbb{Z}$ is integer.

\begin{example}
{\rm Consider the curve $\mathcal{C}_{n,m}$ with affine equation $y^{m}+x^{n}=1$ defined over a finite field $\mathbb{F}_{q}$ of $q$ elements, where $q$ is a power of a prime and $n,m$ are integer numbers greater or equal than 2.

We denote by $F_{n,m}$ the function field $k(x,y)$ of  $\mathcal{C}(n,m)$, where $y^{m}+x^{n}=1$. If $m| q^{2}-1$ then the points $P_{0}=(\alpha,0)$ and $P_{1}=(\beta,0)$ with $\alpha^{m}=1$ and $\beta^{n}=1$ are $\mathbb{F}_{q^{2}}-$ rational points of the curve $\mathcal{C}_{n,m}$  and the root divisors of the elements $x,y \in k(x,y)$ are expressed as $div(y)_{0}=mP_{0}$ and $div(x)_{0}=nP_{1}$. It is a cyclic covering of $\mathbb{P}^{1}(\mathbb{F}_{q^{2}})$ of degree $d$, the greatest common divisor of $n$ and $m$. The Galois group is generated by two elements $g_{1}, g_{2}\in PSL(2,q^{2})$ of orders $n$ and $m$ respectively}.

\end{example}



\begin{teor} \label{classTh}
Fix a genus $g$, a degree $d$ and a partition
$\alpha=(\alpha_{1},\alpha_{2},\ldots, \alpha_{m})$ of $d$. This corresponds to a
branched covering $C_{\alpha,d}$ of $\mathbb{P}^{1}$, with $r=d+m+2(g-1)$, with monodromy above $\infty$ given by $\alpha$, and no other specified simple branch points.
We can classify all such possible degree $d$, branched coverings of $\mathbb{P}^{1}$ 
by a genus $g$ connected Riemann surface by realizing every possible automorphism group $Gal(F/F_{0})$.
\end{teor}
{\it Proof.} 
Every degree $d$ cyclic covering $C_{d}$ of the projective line, after a birational transformation   corresponds to a cyclic extension $k(x,z)$ of the rational function field $k(x)$ of degree $d$, where $z$ satisfies an algebraic equation:
\begin{equation}
z^{d}:=\prod_{i=1}^{m}(x-\rho_{i})^{\alpha_{i}}, \ \ \ \ 0<\alpha_{i}<d.
\end{equation}
If $n:=\sum_{i=1}^{s} \alpha_{i}\equiv 0 (d)$ then the place at $\infty$ does not ramify at the above extension. The only places of $k(x)$ that are ramified are the places $p_{i}$ that correspond to the points $x=\rho_{i}$. 

\noindent If the covering ramifies only at 0 and there is no other branching, then $Gal(k(C_{d}/k(x))\cong \mathbb{Z}_{d}$. In this case there is no ramification over $\infty$ (i.e. $\alpha=(1^{d})$).

\noindent If $C_{d}$  ramifies at $\infty$, we recover all possible
cases by projecting $G/\mathbb{Z}_{d}$  into the known finite
subgroups of $PGL(2,k)$, that constitutes the automorphism group of
the rational function field. If $k$ is algebraically closed, as we are
only interested in enumerating all possible conjugacy classes that can
appear and as the base field $k$ contains all roots of unity, it is enough to determine all finite subgroups of $PSO(2)$. By the classification theorem of finite simple groups of $SO(3)$, these are the ternary groups: 
$\mathbb{Z}_{2}\times \mathbb{Z}_{2}$, $D_{n}$, $A_{4}$,  $A_{5}$ and $S_{4}$. If $k$ is arbitrary, $PGL(2,k)$ is contained in $PGL(2,\overline{k})$ and by the same argument we conclude.

\cqd

\subsection{Triangle curves}
The problem of enumerating branched coverings of $\mathbb{P}^{1}$ is
reduced to the combinatorial problem of studying factorizations
$\sigma=\tau_{1}\ldots \tau_{r}$ into $r$ transpositions for any $d$,
$\sigma$ and $r$. The case in which there is no ramification at
$\infty$ corresponds to the partition $\alpha=(1^{d})$. Hurwitz
numbers enumerate non-singular, genus $g$ curves  expressible as
$d-$sheeted coverings of $\mathbb{P}^{1}$, 
with specified branching above one point, simple branching over other specified points and no other branching.
In this section, we study coverings of $\mathbb{P}^{1}$ with ramification at 3 points. 

\begin{defi}A complex algebraic curve $C$ will be termed triangle curve if it admits a finite group of automorphisms $G<Aut(C)$ so that $C/G\cong \mathbb{P}^{1}$ and the natural projection $$C\rightarrow C/G,$$ ramifies over 3 values, say 0, 1, and $\infty$.
\end{defi}

If the branching orders at these points are $p,q$ and $r$ we will say that $C/G$ is an orbifold of type $(p,q,r)$. Due to a celebrated theorem of Belji, triangle curves are known to be defined over a number field. 


If the number of orbifold points is at least 3, we have the following possibilities for the orders of the orbifold points:
$(2,2,n), $ for some $n\geq 2$, $(2,2,3)$, $(2,3,4)$ and $(2,3,5)$.

The corresponding fundamental group is the dihedral, tetrahedral or icosahedral group repectively, and the universal covering is $\mathbb{P}^{1}$. Any finitely generated discrete subgroup $G$ of $PSL(2,\mathbb{R})$, is the fundamental group of an orbifold and hence it has a presentation of the form:

$$G=<a_{1},b_{1},\ldots, a_{g},b_{g},c_{1},\ldots, c_{k}|\, c_{1}^{n_{1}}=c_{2}^{n_{2}}$$

$$ \ldots=c_{k}^{n_{k}}=1, \ \ [a_{1},b_{1}][a_{2},b_{2}]\cdots[a_{g},b_{g}]c_{1}c_{2}\ldots c_{k}=1>.$$

\begin{prop} 
 All the coverings of $\mathbb{P}^{1}$ that ramify over 3 points are encoded by the partitions of 3 parts: $(2,2,n), $ for some $n\geq 2$, $(2,3,3)$, $(2,3,4)$ and $(2,3,5)$.
\end{prop}

{\it Proof.} All the coverings of $\mathbb{P}^{1}$ that ramify over 3 points are induced by 
the two groups generated by 2 of the 3 transpositions of $S_{3}$, that is  $H_{1}=<(12), (23)>$ and the group $H_{2}=<(23),(13)>$. Thus the  covering whose ramification is given by the 3 permutations $(12)$, $(23)$ and $(12)(23)$ in $H_{1}=<(23), (13)>$  has two simple branch points corresponding to the two transpositions and a branch point with multiplicity at least 3 corresponding to the permutation of order 3,  $(132)$ and all its powers. These coverings have ramification type above $\infty$ defined by the partition $(2,2,n)$  corresponding to the orders of the three orbifold points and the Galois group is the dihedral group $D_{n}$. 
If we consider the  group  $H_{2}=<(23), (13)>$, we recover  the other possible triangle groups $A_{4}$, $S_{4} $, and $A_{5}$ corresponding to the partitions (2,3,3), (2,3,4) and (2,3,5).
\cqd

\section{Enumerative geometry of coverings of $\mathbb{P}^{1}$}\label{sec4}
\subsection{Coverings of $\mathbb{P}^{1}$ with specified ramification above 0 and $\infty$}
Let $d$ and $g\geq 0$ be integer numbers  representing the degree and
the genus of a covering of $\mathbb{P}^{1}$, and let $\lambda$ and
$\rho$ be partitions of $d$ prescribing the profiles of the covering over $0$ and $\infty$. 
Each covering corresponds to a combinatorial object: a labelled graph with $d$ vertices, $d+g-1$ edges and without loops.

 A connected labeled floor diagram $\mathcal{D}$ of degree $d$ and genus $g$ is a connected oriented graph $G=(V,E)$ on linearly ordered $d-$ element vertex $V$, together with a weight function $w:\mathbb{E}\rightarrow \mathbb{Z}_{>0}$ such that the edge set $E$ consists of $d+g-1$ edges, and each edge in $E$ is directed from a vertex $u$ to a vertex $v>u$, expressing compatibility with linear ordering on $V$. The multiplicity $\mu(\mathcal{D})$ is the product of the squares of $w(e)$ for every edge $e\in E$, that is, 
 $$\mu(\mathcal{D})=\prod_{e\in E}(w(e))^{2}.$$

 \begin{prop} Given $\lambda$ and $\rho$ two partitions of $d$, the set of irreducible complex algebraic curves
\begin{itemize}
\item
of degree $d$ and genus $g$ passing through a generic configuration of
$2d-1+g+l(\rho)$ points in $\mathbb{C}^{2}$ 
\item
having tangency to the $x-$ axis for a given collection
$\mathcal{P}_{\lambda}$ of $l(\lambda)$ points in $\mathbb{C}\times
\{0\}$ and other $l(\rho)$ points 
\end{itemize}
coincides with the set of irreducible plane curves $\gamma$ of given degree and genus realizable as $d-$sheeted coverings of $\mathbb{P}^{1}(\mathbb{C})$ with ramification type at 0 and $\infty$ described by the partitions $\lambda$ and $\rho$ and simple ramification over the specified collection of points $\mathcal{P}_{\lambda}$. 
 \end{prop}

 {\it Proof.}
 As we showed in Lemma \ref{lema1}, given an irreducible plane algebraic curve, if we impose the curve to pass through a generic point in the plane, we get a $d-$ sheeted branched covering of $\mathbb{P}^{1}(\mathbb{C})$, by projecting onto the $x-$axis. Furthermore, we can recover the $y-$coordinates by taking $d-$roots of the $x-$coordinate. If the curve has a tangency to the $x-$axis at a generic point of affine coordinates $(x_{i},0)$, the corresponding sheeted covering is branched at this point with the same multiplicity.\cqd
 
\begin{remark}
The authors proved in \cite{FM} that the Gromov-Witten invariant
$N_{d,g}$ representing the number of irreducible curves of degree $d$
and genus $g$ passing through a fixed generic configuration of
$3d+g-1$ points on $\mathbb{P}^{2}$, can be obtained by summing the
product of corresponding multiplicities 
$\mu(\mathcal{D})\cdot\nu(\mathcal{D})$ over all labeled floor diagrams $\mathcal{D}$ of degree $d$ and genus $g$. While the numbers $N_{d,g}(\lambda, \rho)$ count irreducible plane curves $\gamma$ of given degree and genus realizable as $d-$sheeted coverings of $\mathbb{P}^{1}$ with ramification type at 0 and $\infty$ described by the partitions $\mu$ and $\nu$ and simple ramification over other specified points. If $\lambda$ and $\rho$ are two partitions with $|\lambda|+|\rho|=d$, the number $N_{d,g}(\lambda,\rho)$ can be obtained by summing the multiplicities $\mu(\mathcal{D})\nu_{\lambda,\rho}(\mathcal{D})$, where $\nu_{\lambda, \rho}(\mathcal{D})$ is the multiplicity of a certain combinatorial decoration of a labelled floor diagram $\mathcal{D}$.  
\end{remark}


\subsection{Coverings of $\mathbb{P}^{1}$ with 4 or more branch points}
Let $p_{1}\ldots, p_{r}$ be points in $\mathbb{P}^{1}$ and $(s_{1},\ldots, s_{r})$ a set of $r$ permutations defined up to conjugation in $S_{d}$  such that $s_{1}s_{2}\cdot\ldots s_{r}=1$, and the corresponding cycle types are given by partitions $(\eta^{1}\ldots, \eta^{r})$ of $d$, defining the ramification profile  over $p_{i}$.

 There are only finitely many coverings $H^{\mathbb{P}^{1}}_{d}(\eta^{1},\ldots, \eta^{r})$  of the projective line up to isomorphism by smooth connected curves of specified degree and genus, and monodromy $\eta^{i}$ at $p_{i}$. Each covering $\pi$ has a finite group of automorphisms $Aut\,(\pi)$.
This number can be computed by operating in the group algebra  $\mathbb{Q}S_{d}=\{\sum_{\sigma \in S_{d}}\lambda_{\sigma}\sigma,\ \ \lambda_{\sigma}\in \mathbb{Q}\}$ of $S_{d}$.
Let $\mathcal{P}(d)$ denote the set of partitions of $d$ indexing the irreducible representations of $S_{d}$. 
The class algebra $\mathcal{Z}_{d}\subset \mathbb{Q}S_{d}$ is the center of the group algebra.  Let $c_{\eta}\in \mathcal{Z}_{d}$ be the conjugacy class corresponding to the partition $\eta$, then:
\begin{equation}
H^{\mathbb{P}^{1}}_{d}(\eta^{1},\ldots, \eta^{r})=\frac{1}{d!}\left[ C_{(1^{d})}\right]\prod C_{\eta^{i}},
\end{equation}
where $C_{(1^{d})}$ stands for the coefficient of the identity class.

A labelled partition of $d$ is a partition in which the terms are considered distinguished. For example, there are ${7 \choose 3}$ ways of splitting the labelled partition $\alpha=[1^{7}]$ into two labelled partitions $\beta=[1^{3}]$ and $[\gamma]=[1^{4}]$ ($\gamma=\alpha\backslash \beta$).

The $\mathbb{Q}-$algebra structure of $\mathbb{Q}S_{d}$ is given by the unit $u$ and the multiplication $m: \mathbb{Q}S_{d}\otimes \mathbb{Q}S_{d}\rightarrow \mathbb{Q}S_{d}$ defined by the formula: $[\lambda]\otimes [\mu]=\bigoplus_{\rho} c_{\lambda\mu}^{\rho}[\rho]$, where we call the structure constants as $c_{\lambda\mu}^{\rho} \in \mathbb{N}$. 
  If we look at the group algebra $\mathbb{Q}S_{d}$ from a Hopf algebra perspective, an additive basis of $\mathbb{Q}S_{d}$ is indexed by partitions $\{[\lambda]\}_{\lambda \in \mathcal{P}(d)}$. In particular there is an isomorphism with the Hopf algebra of Schur functions and with the Hopf algebra of irreducible representations of $GL(\mathbb{C}^{d})$. Let us call by $k_{\lambda\mu\rho}$ the structure constants for the coproduct $\triangle[\eta]=\sum k^{\eta}_{\lambda, \mu} [\lambda]\otimes [\mu]$ and $S$ the antipode, that is,   $S(\sigma)=\sigma^{-1}$,  $\forall \sigma \in S_{d}$.
 Then the coefficients $k_{\lambda\mu\rho}$ for the coproduct $\triangle[\eta]$ correspond to the structure constants of the dual Hopf algebra $\mathbb{Q}S_{d}$ that are known as Kronecker coefficients.
  

\begin{prop} 

\begin{enumerate}
\item \label{1}The structure constants $c_{\lambda, \mu}^{\eta}$ for the product $m$ of the Hopf algebra $\mathbb{Q}S_{d}$ are the Littlewood-Richardson coefficients.
\item \label{2} The coefficients $k_{\lambda\mu\rho}$ of the coproduct $\triangle$ are the structure constants for the dual Hopf algebra that are known as Kronecker coefficients. 
\end{enumerate}
\end{prop}
{\it Proof.} 
(\ref{1})  In terms of irreducible representations of
$GL(\mathbb{C})$, a partition $\eta$ corresponds to a finite
irreducible representation that we denote as $V(\eta)$. Since
$GL_{d}(\mathbb{C})$ is reductive, any finite dimensional
representation decomposes into a direct sum of irreducible
representations, and the structure constant $c^{\eta}_{\lambda, \mu}$
is the number of times that a given irreducible representation
$V(\eta)$ appears in an irreducible decomposition of
$V(\lambda)\otimes V(\mu)$. These are known as Littlewood-Richardson
coefficients, since they were the first to give a combinatorial
formula encoding these numbers (see \cite{Fu}). 

There is a description of the Littlewood-Richardson coefficients in terms of Young diagrams.
For example, if we consider the partition $\alpha=(5,3,3,1)$, its Young diagram is:
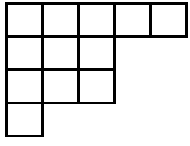
\begin{equation*}
\begin{Young}
 &  &   &   & \cr
 &  & \cr
 &  &  \cr
 \cr
\end{Young}
\end{equation*}

If we represent the partitions $\lambda, \mu, \eta$ by the corresponding Young diagrams, the coefficient $c^{\eta}_{\lambda,\mu}$ represent the number of ways to fill the boxes $\eta\backslash \lambda$, with one integer in each box, so that the following conditions are satisfied:
\begin{itemize}
\item The entries in any row are weakly increasing from left to right.
\item The entries in each column are strictly increasing from top to botton.
\item The integer $i$ occurs exactly $\mu_{i}$ times.
\item For any $p$ with $1\leq p < \sum \mu_{i}$, and any $i$ with $1\leq i <n$, the number of times $i$ occurs in the first $p$ boxes of the ordering is at least as large as the number of times that $i+1$ occurs in these first $p$ boxes.
If we regard an $n-$tuple of parts of the partition $\lambda$ as a point $(\lambda_{1},\ldots, \lambda_{n})$ in $\mathbb{R}^{n}$, then the point corresponding to the partition $\eta$ must be 
in the convex hull of the points $\lambda+ \mu_{\sigma}$, where $\sigma$ varies over the symmetric group $S_{d}$ and $\mu_{\sigma}$ denotes  $\mu_{\sigma}=(\mu_{\sigma(1)},\ldots, \mu_{\sigma(n)})$.
\end{itemize}
(\ref{2}) In terms of the Hopf algebra $\Lambda$ of Schur functions, let $s_{\lambda}$ the Schur function indexed by the partition $\lambda$, we have  $s_{\lambda}\cdot s_{\mu}= \sum_{\nu} c_{\lambda\mu}^{\nu}s_{\nu}$ for the product and we get the coefficients $k_{\lambda \mu}^{\rho}$ as the structure constants of the dual Hopf algebra $\Lambda^{*}$. These are known as Kronecker coefficients, (see \cite{Ma} and \cite{SLL}).


\cqd



\subsection{Connection with the moduli space of curves: Enumeration of Hurwitz numbers}


The symmetric group $S_{n}$ acts on $\mathbb{C}[x_{1},\ldots, x_{n}]$ by 
$$(s\cdot f)(x_{1},\ldots, x_{n})=f(x_{s(1)},\ldots, x_{s(n)}), \ \ {\rm{for}}\ \ s\in S_{n},\ \ f\in \mathbb{C}[x_{1},\ldots, x_{n}].$$
We can view this as the action of $S_{n}$ on $\mathcal{P}(\mathbb{C}^{n})$ arising from the representation of $S_{n}$ on $\mathbb{C}^{n}$ as permutation matrices, with $x=[x_{1},\ldots, x_{n}] \in \mathbb{C}^{n}$.

If $V_{S_{n}}$ is the variety parametrizing curves with Galois group $S_{n}$ then the subvariety of invariants by the action of finite subgroups of $S_{n}$ defines an stratification of the ambient variety $V_{S_{n}}$.

Fix $m$ points $q_{1}\ldots, q_{m}$ in $\mathbb{P}^{1}$ and a conjugacy class $\sigma=(l_{1})\ldots (l_{m})$ in $S_{d}$. Consider the corresponding covering $p:C\rightarrow \mathbf{ P}^{1}$ with ramification type prescribed by the partition $\mu=(l_{1}, \ldots,l_{m})\in \mathcal{P}(d)$ with the integers $l_{i}$ for $i=1, \ldots, m$ ordered by non-decreasing order. 
The pre-image $p^{-1}(\infty)=\sum_{i=1}^{m}l_{i}q_{i}$, defines a divisor on $C$. Let $k(C_{\mu,d})$ be the function field of the curve $C$, we have that $k(C_{\mu,d})\cong k(a_{1},\ldots,a_{m})$, where \begin{equation} \label{eq2} y^{n}=\prod_{i=1}^{m}(x-q_{1})^{l_{1}} \ldots (x-q_{m})^{l_{m}}=\sum_{i=0}^{m}a_{i}x^{i}, \end{equation}
and the coefficients $a_{0},\ldots a_{m}$ are symmetric polynomials of $q_{i}$ multiplied by $(-1)^{s-i}$. The partition $\mu$ gives information on the cycle structure of the permutation $\sigma$.

\begin{remark} Let $\gamma_{k}=\sharp \{l_{i}=k\}$ be the number of
  times the multiplicity corresponding to the integer $k$ is
  realized. If we fix $q_{1},\ldots, q_{m}$ points in
  $\mathbb{P}^{1}$, and we assume the image of $z$ is $p(z)=\infty$, then for each degree $n$,  the number of branched coverings with  the same monodromy type above $\infty$ defined by the partition $\mu$, that is the number of coverings defined by the equation ($\ref{eq2}$),  coincides with the number $\frac{m!}{\prod_{i=1}^{m}\gamma_{i}}$.
\end{remark}

\begin{remark} If we vary one of the  branch points $q_{i}$  of the curve defined by (\ref{eq1}), we obtain a one dimensional Hurwitz space parameterizing  such  coverings.

\end{remark}

\begin{lem} The variety $V_{d,g,\sigma}$ parametrizing coverings with ramification type corresponding to a conjugacy class in $\sigma \in S_{d}$ is a one dimensional subvariety  of the the variety parametrizing degree $d$ and genus $g$ coverings $V_{d,g}$.
\end{lem}
{\it Proof.} Consider the natural identification of $\sigma_{d}$ with an  element $A_{\sigma}$ in $GL_{d}(k)$ (respectively $SL_{d}(k)$), via a linear representation. This element 
determines an automorphism of the function field given by multiplication of the corresponding matrix representation $A_{\sigma}$ in $GL_{d}(k)$ with the vector field of coordinates $(a_{1},\ldots, a_{m})$. The invariant field $k(a_{1},\ldots, a_{m})^{\sigma_{d}}$ is the quotient surface $C_{\mu,d}/G$ by the group $G$ generated by the corresponding element $A_{\sigma}$ in $GL_{d}(k)$, thus a one dimensional scheme in $V_{d,g}$. \cqd

Let $H_{g,\mu}$ be the Hurwitz number, that is, the number of genus $g$ degree $d$ coverings of $\mathbb{P}^{1}$ with profile $\mu$ over $\infty$ and simple ramification over a fixed set of finite points. The Hurwitz numbers are naturally expressed in terms of tautological intersections in the moduli space ${M}_{g,n}$ of  projective nonsingular curves of genus $g$ and $n$ marked points, and its compactification $\overline{M}_{g,n}$, whose points correspond to projective, connected, nodal curves of arithmetic genus $g$, satisfying a stability condition (due to Deligne and Mumford), and with orbifold singularities if regarded as ordinary coarse moduli spaces. These moduli spaces are irreducible varieties of dimension $3g-3+n$ if $g\geq 2$, smooth if regarded as (Deligne-Mumford) stacks, and with orbifold singularities if regarded as ordinary coarse moduli spaces. 
The Deligne-Mumford compactification $\overline{M}_{g,n}$ of the moduli space of curves comes equipped with a  well-defined Chow intersection ring much like the cohomology ring of a compact manifold.

Much of what we know about $\overline{M}_{g,n}$ comes from intersection
numbers of the tautological $\psi_{i}$, $\kappa$, and $\lambda_{j}$ classes in the tautological subring $R¥(\overline{M}_{g,n})$.
Kontsevich's proof of Witten's conjecture essentially provided a recursive formula for all
intersections of $\psi_{i}$ classes (thus also all $\kappa$ classes), yet the study of relations involving $\lambda_{j}$
classes, or Hodge integrals, is still an active field now closely related to Gromov-Witten
invariants and the combinatorics of Hurwitz numbers.

The intersection theory of $\overline{M}_{g,n}$ must be studied in the orbifold category  or the category of Deligne Mumford-stacks to correctly handle the automorphisms group of the pointed curves. For each marking $i$, there exists a canonical line bundle $\mathbb{L}_{i}$. The fiber  at the stable pointed curve $(C,x_{1},\ldots, x_{n})$ is the cotangent space $T^{*}_{C}(x_{i})$ of $C$ at $x_{i}$. $\mathbb{L}_{i}$ determines a $\mathbb{Q}-$divisor on the coarse moduli space. Let $\psi_{i}$ denote the first Chern class of $\mathbb{L}_{i}$. Witten's conjecture concerns the complete set of evaluations of intersections of the $\psi$ classes:
\begin{equation}\label{eq3} \int_{\overline{M}_{g,n}}\psi_{1}^{k_{1}}\ldots \psi_{n}^{k_{n}}. \end{equation}
The symmetric group $S_{n}$ acts naturally on $\overline{M}_{g,n}$ by permuting the markings. Since the $\psi$ classes are permuted by this $S_{n}$ action, the integral is unchange by a permutation of the exponents $k_{i}$. A notation for these intersections which exploits the $S_{n}$ symmetry is given by:
\begin{equation}\label{eq4} <\tau_{k_{1}}\ldots, \tau_{k_{n}}>_{g}=\int_{\overline{M}_{g,n}}\psi_{1}^{k_{1}}\ldots \psi_{n}^{k_{n}}. \end{equation}

Let the Hodge bundle 
$$\mathbb{E}\rightarrow \overline{M}_{g,n}$$ be the rank $g$ vector bundle with fiber $H^{0}(C,w_{C})$ over the moduli point $(C,p_{1}\ldots, p_{n})$. The $\lambda$ classes are the Chern classes of the Hodge bundle:
$$\lambda_{i}=c_{i}(E)\in H^{2i}(\overline{M}_{g,n},\mathbb{Q}).$$
The $\psi$ and $\lambda$ classes are {\it tautological} classes on the moduli space of curves.The Hurwitz numbers are naturally expressed in terms of tautological intersections in $M_{g,n}$ (see Theorem 2 of \cite{OP}).



Let $H^{g}_{d}$ be the number of such branched coverings that are connected, then the following formula due to Ekedahl, Lando, Shapiro and Vainshtein, expresses Hurwitz numbers in terms of Hodge integrals (see \cite{VGJ}).
$$H^{g}_{\alpha}=\frac{r!}{\sharp Aut(\alpha)}\prod_{i=1}^{m}\frac{\alpha_{i}^{\alpha_{i}}}{\alpha_{i}}\int_{\overline{M}_{g,m}}\frac{1-\lambda_{1}+\ldots,\pm \lambda_{g}}{\prod(1-\alpha_{i}\psi_{i})}$$ 




\begin{example} In the case of genus 0, the formula reads: 
$$H^{0}_{d}=\frac{(2d-2)!}{d!} d^{d-3}.$$

\end{example}

\subsubsection{Connection with the moduli space of curves: counting coverings of $\mathbb{P}^{1}$ over finite fields}

Let $\overline{M}_{g,n}(\mathbb{F}_{p})$ be the moduli space of stable curves of genus $g$ with $n$ marked points defined over the finite field $\mathbb{F}_{p}$ of $p$ elements.

\begin{prop} The number of genus $g$ curves expressible as $d-$sheeted coverings of $\mathbb{P}^{1}$ coincides with the cardinality of $M_{g,m}$ over $\mathbb{F}_{p}$ (up to $\mathbb{F}_{p}$ isomorphism), where $g=\frac{(d-1)(m-2)}{2}$, weighted by the factor  $1/\sharp Aut_{\mathbb{F}_{p}}(C)$.
\end{prop}

{\it Proof.}  Let $C$ be a complete, connected non singular curve with $m$ marked points $p_{1},\ldots, p_{m}$. We obtain a morphism $f:C\rightarrow \mathbb{P}^{1}$ from the linear series attached to the divisor $p_{1}+\ldots+p_{m}$. The branched covering $f$ expresses $C$ in the form $y^{d}=\prod_{i=1}^{m}(x-p_{i})^{l_{i}}$, with profile defined by the partition $(l_{1},\ldots, l_{m})$ of $d$, expressing the monodromy above $\infty$. By Riemman-Hurwitz formula we can compute $m$, the number of different branch points. 
Now the number of polynomials of degree $n=d+m+2\,(g-1)$ with $m$ different roots is the falling factorial polynomial $(p)_{m+1}:=p\,(p-1)\,(p-2)\ldots (p-m)$, divided by the order of the affine transformation group of $\mathbb{A}^{1}=\mathbb{P}^{1}\backslash \infty$, that is, $p^{2}-p$.
\cqd

\subsubsection*{Acknowledgments}
We would like to thank Fei Xu and Vivek Mallick for interesting discussions while the preparation of the work and Joachim Kock for reading the paper and some useful comments. This work has been partially supported by the project MTM2009-10359 \lq\lq M\'etodos combinatorios en Geometr\'ia Aritm\'etica y Geometr\'ia Algebraica".

\end{document}